\documentclass[12pt]{article}
\usepackage{e-jc}

\usepackage[usenames,dvipsnames,table]{xcolor}
\usepackage{graphicx}
\DeclareGraphicsExtensions{.pdf}
\graphicspath{{figs/}} 

\usepackage{amsmath}

\usepackage{cite}
\usepackage[colorlinks=true,citecolor=black,linkcolor=black,urlcolor=blue]{hyperref}
\usepackage[capitalise,noabbrev,nameinlink]{cleveref}
\crefname{ineq}{Inequality}{Inequalities}
\crefname{lemma}{Lemma}{Lemmata}
\creflabelformat{ineq}{#2{\upshape #1}#3} 

\usepackage{enumitem}

\usepackage{mathtools}
\DeclarePairedDelimiter\abs{\lvert}{\rvert}

\newcommand{\cC}{{\mathcal{C}}}
\newcommand{\cE}{{\mathcal{E}}}

\newcommand{\taxa}{\ensuremath{X}}
\def\root{\ensuremath{\rho}}

\usepackage{xspace}
\newcommand{\rSPR}{\textup{rSPR}\xspace} 
\newcommand{\SNPR}{\textup{SNPR}\xspace}
\newcommand{\SNPRZ}{\textup{SNPR$^0$}\xspace}
\newcommand{\SNPRP}{\textup{SNPR$^+$}\xspace}
\newcommand{\SNPRM}{\textup{SNPR$^-$}\xspace}
\DeclareMathOperator{\dSNPR}{\ensuremath{d_{\SNPR}}}
\DeclareMathOperator{\drSPR}{\ensuremath{d_{\rSPR}}} 

\def\class{\ensuremath{\mathcal{C}_n}}
\def\trees{\ensuremath{\mathcal{T}_n}}
\def\nets{\ensuremath{\mathcal{N}_n}}
\def\tchinets{\ensuremath{\mathcal{TC}_n}}
\def\tbasednets{\ensuremath{\mathcal{TB}_n}}
\def\retvisnets{\ensuremath{\mathcal{RV}_n}}

\title{On the Subnet Prune and Regraft Distance}
\author{Jonathan Klawitter$^*$ \qquad Simone Linz\thanks{Supported by the New Zealand Marsden Fund.}\\
	\small School of Computer Science\\[-0.8ex]
	\small University of Auckland\\[-0.8ex] 
	\small Auckland, New Zealand\\
	\small\tt jo.klawitter@gmail.com, s.linz@auckland.ac.nz\\
}

\dateline{May 11, 2018}{Mar 13, 2019}{Apr 5, 2019}
\Copyright{The authors. Released under the CC BY license (International 4.0).}

\begin{document} 

\maketitle

\begin{abstract} 
\pdfbookmark[1]{Abstract}{Abstract}
Phylogenetic networks are rooted directed acyclic graphs that represent evolutionary relationships
between species whose past includes reticulation events such as hybridisation and horizontal gene
transfer. 
To search the space of phylogenetic networks, the popular tree rearrangement operation rooted
subtree prune and regraft (rSPR) was recently generalised to phylogenetic networks.
This new operation -- called subnet prune and regraft (SNPR) -- induces a metric on the space of all
phylogenetic networks as well as on several widely-used network classes.
In this paper, we investigate several problems that arise in the context of computing the
SNPR-distance.
For a phylogenetic tree $T$ and a phylogenetic network $N$, we show how this distance can be
computed by considering the set of trees that are embedded in $N$ and then use this result to
characterise the SNPR-distance between $T$ and $N$ in terms of agreement forests.
Furthermore, we analyse properties of shortest SNPR-sequences between two phylogenetic networks $N$
and $N'$, and answer the question whether or not any of the classes of tree-child,
reticulation-visible, or tree-based networks isometrically embeds into the class of all phylogenetic
networks under SNPR.
\end{abstract} 

\section{Introduction}
\label{sec:introduction}

Many algorithms that have been developed to reconstruct phylogenetic trees from molecular sequence
data require a (heuristic) search of the space of all phylogenetic trees~\cite{Fel04}. 
To this end, local rearrangement operations, such as nearest neighbor interchange, subtree prune and regraft, 
and tree bisection and reconnection, have been introduced that induce metrics on the space of phylogenetic trees~\cite{SOW96}.
More recently, rooted phylogenetic networks, which are leaf-labelled rooted directed acyclic graphs,
have become increasingly popular in the analysis of ancestral relationships between species whose
past includes speciation as well as reticulation events such as hybridisation and horizontal gene transfer~\cite{Gus14,HRS10}.
In particular, each vertex in a rooted phylogenetic network whose in-degree is at least two
represents a reticulation event and is referred to as a \emph{reticulation}. 
In comparison to tree space, the space of phylogenetic networks is significantly larger and
searching this space remains poorly understood although the above-mentioned rearrangement
operations on phylogenetic trees have been generalised to rooted (and unrooted) phylogenetic
networks~\cite{BLS17,FHMW17,GvIJLPS17,HLMW16,HMW16,JJEvIS17,Kla17}.
 
The goal of this paper is to advance our understanding of the subnet prune and regraft (SNPR)
operation~\cite{BLS17} with a particular focus on the induced distance. 
For two phylogenetic networks, this distance equates to the minimum number of SNPR operations that
are required to transform one network into the other one.
SNPR generalises the rooted subtree prune and regraft (rSPR) operation~\cite{AS01,BS05,SOW96} from
rooted phylogenetic trees to rooted phylogenetic networks. 
A second generalisation of the rSPR operation from trees to networks was recently introduced by Gambette et al.~\cite{GvIJLPS17}.
Both generalisations are similar in the sense that they allow horizontal as well as vertical rearrangement moves. 
From a practical perspective, the space of phylogenetic networks can be searched horizontally in tiers, 
where a tier contains all phylogenetic networks with a fixed number of reticulations, 
as well as vertically among different tiers since a single operation can increase or decrease the number of reticulations by at most one.
On the other hand, there are also subtle differences between the two operations. 
While SNPR is defined on rooted phylogenetic networks that allow for parallel edges~\cite{BLS17},
the generalisation of rSPR to networks as introduced by Gambette et al.~\cite{GvIJLPS17} is defined on networks that do not allow for parallel edges. 
Moreover, the latter operation allows for the switching of a parent vertex (referred to as a \emph{tail moves}) 
and for the switching of a child of a reticulation (referred to as a \emph{head moves}) while SNPR only allows for tail moves. 
Under SNPR, tail moves are sufficient to establish that the operation induces a metric on the space of all rooted phylogenetic networks. 
Moreover, SNPR also induces a metric on the space of several popular classes of phylogenetic networks, such as tree-child, reticulation-visible, and tree-based networks~\cite{CRV09,FS15},
regardless of whether or not one restricts to subclasses of these networks that have a fixed number of reticulations.

Since computing the rSPR-distance between two phylogenetic trees is NP-hard~\cite{BS05}, 
it is not surprising that calculating the SNPR-distance as well as the distance induced by the operation introduced by Gambette et al.~\cite{GvIJLPS17} 
and further investigated by Janssen et al.~\cite{JJEvIS17} is also NP-hard. 
In this paper, we investigate  problems that arise in the context of computing the SNPR-distance.
Bordewich et al.~\cite{BLS17} established several bounds on the SNPR-distance and showed that, 
for a rooted phylogenetic tree $T$ and a rooted phylogenetic network $N$, 
the SNPR-distance $\dSNPR(T, N)$ between $T$ and $N$ is equal to the number of reticulations in $N$ if $T$ is embedded in $N$. 
In the first part of this paper, we extend their result by showing how $\dSNPR(T, N)$ can be computed regardless of whether or not $T$ is embedded in $N$. 
Roughly speaking, the problem of computing the SNPR-distance is equivalent to computing the minimum rSPR-distance between all tree pairs consisting of $T$ and a tree embedded in $N$. 
Hence, one way of computing $\dSNPR(T, N)$ is by repeatedly solving the rSPR-distance problem between two trees.
We use this result to show that computing $\dSNPR(T, N)$ is fixed-parameter tractable.
We then show that $\dSNPR(T, N)$ can also be characterised in terms of agreement forests. 
The notion of agreement forests is the underpinning concepts for almost all theoretical results as well as 
practical algorithms that are related to computing the rSPR-distance between two rooted phylogenetic trees~\cite{BS05,CFS15,WBZ16,Wu09}.
We extend this notion to computing $\dSNPR(T, N)$, which allows us to work directly on $T$ and $N$ instead of different tree pairs. 
In the second part of this paper, we turn to problems that are related to finding shortest
SNPR-sequences for two rooted phylogenetic networks $N$ and $N'$ with $r$ and $r'$ reticulations, respectively.
In particular, we are interested in the properties of networks that a shortest SNPR-sequence from $N$ to $N'$ contains besides $N$ and $N'$. 
For example, if there is always a sequence with the property that each network in the sequence has 
at least $\min(r, r')$ and at most $\max(r, r')$ reticulations, then this might have positive
implications in devising practical search algorithms because the search space could be pruned appropriately.
Surprisingly, we find that, even if $r = r'$, it is possible that every shortest SNPR-sequence for
$N$ and $N'$ contains a network with strictly more than $r'$ reticulations. 
Moreover, for each $r$ with $r \geq 1$, there exist two rooted phylogenetic networks that both have $r$ reticulations and
for which every shortest SNPR-sequence contains a rooted phylogenetic tree. 

The paper is organised as follows. The next section contains notation and terminology that is used throughout the rest of this paper. 
\Cref{sec:treeNetwork} establishes a new result that equates the SNPR-distance between a phylogenetic tree $T$ 
and a phylogenetic network $N$ to the rSPR-distance between pairs of trees.
This result is used in \cref{sec:maf} to characterise the SNPR-distance between $T$ and $N$ in terms of agreement forests. 
We then investigate properties of shortest SNPR-sequences between two phylogenetic networks in \cref{sec:shortestPaths}.
We end this paper with some concluding remarks in \cref{sec:discussion}.

\section{Preliminaries}
\label{sec:preliminaries}

This section provides notation and terminology that is used in the remainder of the paper. 
In particular, we will introduce notation in the context of phylogenetic networks as well as the
SNPR operation. Throughout this paper, $X = \{1, 2, \ldots, n\}$ denotes a finite set.

\paragraph{Phylogenetic networks.}
\pdfbookmark[2]{Phylogenetic networks}{PhyNets}
A \emph{rooted binary phylogenetic network} $N$ on $X$ is a rooted directed acyclic graph with the following vertices:
\begin{itemize}
  \item the unique \textit{root} $\rho$ with in-degree zero and out-degree one,
  \item \emph{leaves} with in-degree one and out-degree zero bijectively labelled with $X$,
  \item \emph{inner tree vertices} with in-degree one and out-degree two, and
  \item \emph{reticulations} with in-degree two and out-degree one. 
\end{itemize}
The \emph{tree vertices} of $N$ are the union of the inner tree vertices, the leaves and the root. 
An edge $e = (u, v)$ is called \emph{reticulation edge}, if $v$ is a reticulation, and
\emph{tree edge}, if $v$ is a tree vertex.
The set $X$ is referred to as the \emph{label set} of $N$ and is sometimes denoted by $L(N)$. 
Following Bordewich et al.~\cite{BLS17}, we allow edges in $N$ to be in \emph{parallel}, that is,
two distinct edges join the same pair of vertices. 
Also note that our definition of the root is known as \emph{pendant root} \cite{BLS17} and it
differs from another common definition where the root has out-degree two.
Our variation serves both elegance and technical reasons.

Let $N$ be a rooted binary phylogenetic network on $X$. 
For two vertices $u$ and $v$ in $N$, we say that $u$ is a \emph{parent} of
$v$ and $v$ is a \emph{child} of $u$ if there is an edge $(u, v)$ in $N$.
Similarly, we say that $u$ is \emph{ancestor} of $v$ and $v$ is \emph{descendant} of $u$ if there is a directed path from $u$ to $v$ in $N$. 
The vertices $u$ and $v$ are \emph{siblings} if they have a common parent.
Lastly, if $u$ and $v$ are siblings and also leaves, we say they form a \emph{cherry}. 

A \emph{rooted binary phylogenetic tree} on $X$ is a rooted binary phylogenetic network that has no reticulations. 

To ease reading, we refer to a rooted binary phylogenetic network (resp. rooted binary phylogenetic tree)
on $X$ simply as a phylogenetic network or network (resp. phylogenetic tree or tree).
Furthermore, let $\nets$ denote the set of all phylogenetic networks on $X$ and let $\trees$ denote
the set of all phylogenetic trees on $X$ where $n = \abs{ X }$.

Let $G$ be a directed graph.
A \emph{subdivision} of $G$ is a graph that can be obtained from $G$ by subdividing each edge of $G$ with zero or more vertices.
Let $N \in \nets$.
We say $G$ has an \emph{embedding} into $N$ if there exists a subdivision of $G$ that is a subgraph of $N$.
Note that such an embedding maps a labelled vertex of $G$ to a vertex of $N$ with the same label.
Furthermore, we say an embedding of $G$ into $N$ \emph{covers} a vertex $v$ (resp. an edge $e$) of $N$ if a vertex (resp. an edge) of the subdivision of $G$
is mapped to $v$ (resp. $e$) by the embedding.

Let $T \in \trees$ and $N \in \nets$. We say $N$ \emph{displays} $T$ if $T$ has an embedding into $N$.
The set of all phylogenetic trees that are displayed by $N$ is denoted by $D(N)$.

\paragraph{Classes of phylogenetic networks.}
\pdfbookmark[2]{Network classes}{NetworkClasses} 
Let $N \in \nets$.
The network $N$ is a \emph{tree-child} network if each of its non-leaf vertices has a tree vertex as child. 
A vertex $v$ of $N$ is called \emph{visible} if there is a leaf $l$ in $N$ such that every directed path from the root of $N$ to $l$ traverses $v$.
We say that $N$ is a \emph{reticulation-visible} network if every reticulation of $N$ is visible.
Lastly, $N$ is \emph{tree based} if there exists an embedding of a phylogenetic tree $T \in \trees$ into $N$ that covers every vertex of $N$.
For a fixed $n$, the class of tree-child networks is denoted by $\tchinets$, 
of reticulation-visible networks by $\retvisnets$, and of tree-based networks by $\tbasednets$.
Each tree-child network is also a reticulation-visible network~\cite{HRS10} and
each reticulation-visible network is also a tree-based network~\cite{GGLVZ15,FS15}.

\paragraph{SNPR.} 
\pdfbookmark[2]{SNPR}{SNPR}
Let $N \in \nets$ with root $\root$ and let $e = (u, v)$ be an edge of $N$. 
Bordewich et~al.~\cite{BLS17} introduced the \emph{SubNet Prune and Regraft (\SNPR)} operation that
transforms $N$ into a phylogenetic network $N'$ in one of the following three ways:

\begin{enumerate}[leftmargin=*,label=(SNPR$-$)]
  \item[(\SNPRZ)] If $u$ is a tree vertex (and $u \neq \root$), then delete $e$, suppress $u$,
  subdivide an edge that is not a descendant of $v$ with a new vertex $u'$, and add the edge $(u',
  v)$.
  \item[(\SNPRP)] Subdivide $(u, v)$ with a new vertex $v'$, subdivide an edge in the resulting
  network that is not a descendant of $v'$ with a new vertex $u'$, and add the edge $(u', v')$.
  \item[(\SNPRM)] If $u$ is a tree vertex and $v$ is a reticulation, then delete $e$, and suppress
  $u$ and $v$.
\end{enumerate}

\begin{figure}[htb]
 \centering
 \includegraphics{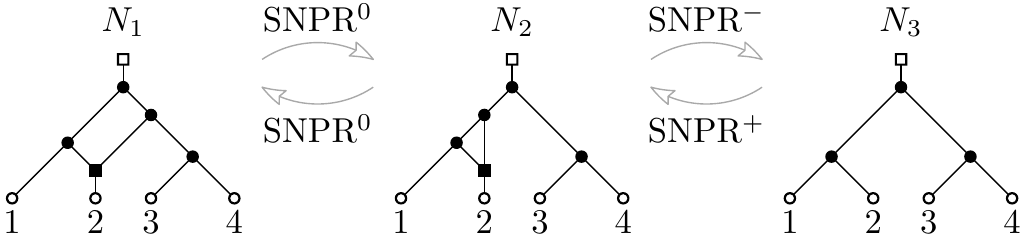}
 \caption{The phylogenetic network $N_2$ can be obtained from $N_1$ by an \SNPRZ and
 the phylogenetic network $N_3$ can be obtained from $N_2$ by an \SNPRM.
 Both operations have a corresponding \SNPRZ and \SNPRP, respectively, that 
 reverses the transformation.}
 \label{fig:SNPR:example}
\end{figure}

In what follows, we sometimes need to specify which of the three operations we consider, in which
case we use $0$, $+$, or $-$ as a superscript to indicate the type of operation.
The three types of operations are illustrated in \cref{fig:SNPR:example}.
Note that an \SNPRZ does not change the number of reticulations, while an \SNPRM decreases
it by one and an \SNPRP increases it by one. 
Lastly, it is worth noting that the well known \rSPR operation~\cite{BS05} on phylogenetic trees is
a restriction of \SNPR in which $N$ and $N'$ are phylogenetic trees and $N$ is transformed into $N'$
by \SNPRZ operations.

\paragraph{SNPR-distance.}
\pdfbookmark[2]{SNPR-distance}{distance} 
Let $N, N' \in \nets$. 
An \emph{\SNPR-sequence} from $N$ to $N'$ is a sequence 
	$$\sigma =(N = N_0, N_1, N_2, \ldots, N_k = N')$$ 
of phylogenetic networks such that, for each $i \in \{ 1, 2, \ldots, k \}$, we can obtain
$N_i$ from $N_{i-1}$ by a single SNPR. The \emph{length} of $\sigma$ is $k$.

Now, let $\cC$ be a class of phylogenetic networks. 
Then $\cC$ is said to be \emph{connected} (under SNPR) if, for all pairs $N$ and $N'$ of
networks in $\cC$, there exists an SNPR-sequence $\sigma$ from $N$ to $N'$ and each network in $\sigma$ is in $\cC$.
Moreover, if $\cC$ is connected, then the \emph{\SNPR-distance} between two elements in $\cC$, say $N$ and $N'$, 
is the length of a shortest SNPR-sequence from $N$ to $N'$ with the property that each network of the sequence is in $\cC$.
This distance is denoted by $d_{\SNPR_{\cC}}(N, N')$ or more simply by $\dSNPR(N, N')$ if the class under
consideration is clear from the context. 
Finally, let $\cC$ and $\cC'$ be two connected classes of phylogenetic networks such that all elements
in $\cC$ are also contained in $\cC'$.
We say that $\cC$ \emph{isometrically embeds} into $\cC'$ if 
$d_{\SNPR_{\cC}}(N,N') = d_{\SNPR_{\cC'}}(N, N')$ for all pairs $N$ and $N'$ of networks
in $\cC$.

For the SNPR-distance to be a metric on a class of networks, the class has to be connected under SNPR and the
SNPR operation has to be reversible, that is, if a phylogenetic network $N'$ can be obtained from
a phylogenetic network $N$ by a single SNPR operation, then $N$ can also be obtained from $N'$ by a single SNPR operation.
Bordewich et~al.~\cite{BLS17} established a metric result for the following four classes of phylogenetic networks.
\begin{proposition}[{\hspace{1sp}{\cite[Corollary 3.3]{BLS17}}}]\label{clm:SNPR:metric}
The \SNPR operation induces a metric on each of the classes $\nets$, $\tchinets$, $\retvisnets$, and $\tbasednets$.
\end{proposition}

\section{Characterising the SNPR-distance between a network and a tree}
\label{sec:treeNetwork}
In this section, we characterise the SNPR-distance $\dSNPR(T, N)$ between a phylogenetic network $N$ 
and a phylogenetic tree $T$ in terms of $D(N)$, the set of phylogenetic trees that are displayed by $N$. 
Bordewich et~al.~\cite{BLS17} have shown how to compute this distance if $T$ is displayed by $N$. 
To give a full characterisation of $\dSNPR(T, N)$ regardless of whether or not $T$ is displayed by $N$, we make use of the following three lemmata.

\begin{lemma}[{\hspace{1sp}\cite[Lemma 7.4]{BLS17}}]\label{clm:SNPR:shortestPath:displayedTree} 
Let $N \in \nets$ with $r$ reticulations. Let $T \in D(N)$.
Then 
	$$\dSNPR(T, N) = r\text{.}$$
\end{lemma} 
%

\begin{lemma}[{\hspace{1sp}{\cite[Proposition 7.1]{BLS17}}}]\label{clm:SNPR:isometric:trees}
Let $T, T' \in \trees$. Then
	$$\drSPR(T, T') = \dSNPR(T, T')\text{.}$$ 
Moreover, the class of all phylogenetic trees $\trees$
isometrically embeds into the class of all phylogenetic networks $\nets$ under the SNPR-distance.
\end{lemma} 
%
 
\begin{lemma}[{\hspace{1sp}\cite[Proposition
7.7]{BLS17}}]\label{clm:SNPR:shortestPath:displayedTrees:step}%
Let $N, N' \in \nets$ such that $\dSNPR(N, N') = k$. Let $T \in D(N)$.
Then there exists a phylogenetic tree $T' \in D(N)$ such that $\dSNPR(T, T')\leq k\text{.}$
\end{lemma}

By setting one of the two networks in the previous lemma to be a phylogenetic tree and noting
that the roles of $N$ and $N'$ are interchangeable, the next two corollaries are immediate
consequences of \cref{clm:SNPR:shortestPath:displayedTree,clm:SNPR:shortestPath:displayedTrees:step}.

\begin{corollary}\label{clm:SNPR:shortestPath:distanceToDisplayed}
Let $T \in \trees$ and $N \in \nets$ with $\dSNPR(T, N) = k$.\\
Then $\dSNPR(T, T') \leq k$ for each $T' \in D(N)$.
\end{corollary}
%

\begin{corollary}\label{clm:SNPR:shortestPath:displayedTrees:distance}
Let $N \in \nets$ with $r$ reticulations. Let $T, T' \in D(N)$.\\ Then $\dSNPR(T, T') \leq r$.
\end{corollary}

The main result of this section is the following theorem that characterises the
SNPR-distance between a phylogenetic tree and a phylogenetic network.

\begin{theorem}\label{clm:SNPR:shortestPath:nondisplayedTree}
Let $T \in \trees$. Let $N \in \nets$ with $r$ reticulations. 
Then
	$$\dSNPR(T, N) = \min_{T' \in D(N)}\dSNPR(T, T') + r \text{.}$$ 
  \begin{proof}
	Let $T^*\in D(N)$ such that $\dSNPR(T, T^*)\leq \dSNPR(T, T')$ for each $T'\in D(N)$. 
	Then, by \cref{clm:SNPR:shortestPath:displayedTree,clm:SNPR:isometric:trees}, it follows that
	\begin{equation}\label[ineq]{eq:inequality-one}
	\dSNPR(T, N) \leq d_{\SNPR}(T, T^*) + d_{\SNPR}(T^*, N)
		= \min_{T' \in D(N)}\dSNPR(T, T') + r \text{.}
	\end{equation}
	We next show that 
		$$\dSNPR(T, N) \geq \min_{T' \in D(N)}\dSNPR(T, T') + r\text{.}$$ 
	Suppose that $\dSNPR(T, N) = k$. 
	Let $\sigma = (T = N_0, N_1, N_2,\ldots, N_k = N)$ be an SNPR-sequence from $T$ to $N$. 
	For each $i \in \{1,2,\ldots,k\}$, consider the two networks $N_{i-1}$ and $N_i$ in $\sigma$. 
	If $N_i$ has been obtained from $N_{i-1}$ by applying an SNPR$^+$ operation, then
	$D(N_{i-1})\subseteq D(N_i)$.
	Furthermore, regardless of the SNPR operation used to obtain $N_i$ from $N_{i-1}$ \cref{clm:SNPR:shortestPath:displayedTrees:step} implies that, for each tree
	$T_{i-1} \in D(N_{i-1})$, there exists a tree $T_i$ in $D(N_i)$ such that 
	$\dSNPR(T_{i-1}, T_i)\leq 1$. 
	It is now straightforward to check that we can construct a sequence
	$S = (T_0, T_1, T_2, \ldots, T_{k})$ of phylogenetic trees on $X$ from $\sigma$ that satisfies
	the following properties.
	\begin{enumerate}[label=(\roman*)]  
	  \item For each $i \in \{0, 1, \ldots, k\}$, we have $T_i\in D(N_i)$.
	  \item For each $i \in \{1, 2, \ldots, k\}$, if $N_i$ has been obtained from $N_{i-1}$ by applying
	  an SNPR$^+$ operation, then $T_i = T_{i-1}$.
	  \item For each $i \in \{1, 2, \ldots, k\}$, we have $d_{\SNPR}(T_{i-1}, T_i) \leq 1$.
	\end{enumerate}
	By construction and since $\sigma$ contains at least $r$ \SNPRP operations, there exists a
	subsequence of $S$ of length $k-r$ that is an SNPR-sequence from $T_0$ to $T_k$. 
	Hence, we have $\dSNPR(T, T_k) \leq k-r$. Moreover, noting that $T_k \in D(N)$ it follows from
	\cref{clm:SNPR:shortestPath:displayedTree} that  $\dSNPR(T_k, N) = r$ and, thus,
	\begin{equation}\label[ineq]{eq:inequality-two}
	\begin{split}
	\min_{T'\in D(N)}\dSNPR(T, T')+r & \leq  \dSNPR(T, T_k)+\dSNPR(T_k, N)\\
	&= k - r + r = k = \dSNPR(T, N)\text{.}
	\end{split}
	\end{equation}
	
	Combining \cref{eq:inequality-one,eq:inequality-two} establishes the theorem.
  \end{proof}
\end{theorem}

Given \cref{clm:SNPR:isometric:trees} and \cref{clm:SNPR:shortestPath:nondisplayedTree} and that $\dSNPR(T, T') = \drSPR(T, T')$, 
it is worth noting that the problem of computing the SNPR-distance between a phylogenetic network
and a phylogenetic tree can be reduced to computing the rSPR-distance between pairs of trees. 
Calculating the rSPR-distance between two phylogenetic trees is a well understood problem and
several exact algorithms exist (e.g.~\cite{BS05,WBZ16}).
Furthermore, this problem is known to be fixed-parameter tractable with the rSPR-distance itself as parameter~\cite[Theorem 3.4]{BS05}. 
This means that there exists an algorithm to compute $k = \drSPR(T, T') = \dSNPR(T, T')$ in $f(k)p(n)$ time
where $f$ is a computable function that only depends on $k$ and $p$ is a polynomial function. 
Note that replacing $k$ by a function $f'(k)$ or calling such an algorithm as a black-box at most $f'(k)$ times, yields
again a fixed-parameter tractable algorithm in $k$. 
We use this observation to establish the following theorem.

\begin{theorem}\label{clm:SNPR:treeNetwork:fpt}
Let $T \in \trees$ and $N \in \nets$. Then computing $\dSNPR(T, N)$ is fixed-parameter
tractable when parameterised by $\dSNPR(T, N)$.
  \begin{proof}
    Let $d = \dSNPR(T, N)$ and let $r$ be the number of reticulations of $N$.
    By \cref{clm:SNPR:shortestPath:distanceToDisplayed} we know that $k = \drSPR(T, T') = \dSNPR(T, T') \leq d$ for all $T' \in D(N)$.
	From the observation before the theorem, it follows that computing $\drSPR(T, T')$ is also fixed-parameter tractable when parameterised by $d$.
	Next, note that $\abs{D(N)} \leq 2^r \leq 2^d$, since we know by \cref{clm:SNPR:shortestPath:nondisplayedTree} that $r \leq d$.
	Again, by the observation above, computing $\drSPR(T, T')$ for at most $2^d$ trees $T' \in D(N)$ is still fixed-parameter tractable when parameterised by $d$.
    By \cref{clm:SNPR:shortestPath:nondisplayedTree} $\dSNPR(T, N)$ can be computed by computing $\drSPR(T, T')$ for each $T' \in D(N)$.
    Taken together, this implies that computing $\dSNPR(T, N)$ is fixed-parameter tractable.
  \end{proof}
\end{theorem}

\section{Using agreement forests to characterise the SNPR-distance}
\label{sec:maf}
We now show how agreement forests can be used to characterise the SNPR-distance
between a phylogenetic tree $T$ and a phylogenetic network $N$.
Importantly, this characterisation allows us to compute the SNPR-distance between $T$ and $N$ directly 
without having to compute the rSPR-distance between $T$ and each tree that is displayed by $N$ as suggested by \cref{clm:SNPR:shortestPath:nondisplayedTree}.

\begin{figure}[htb]
  \centering
  \includegraphics{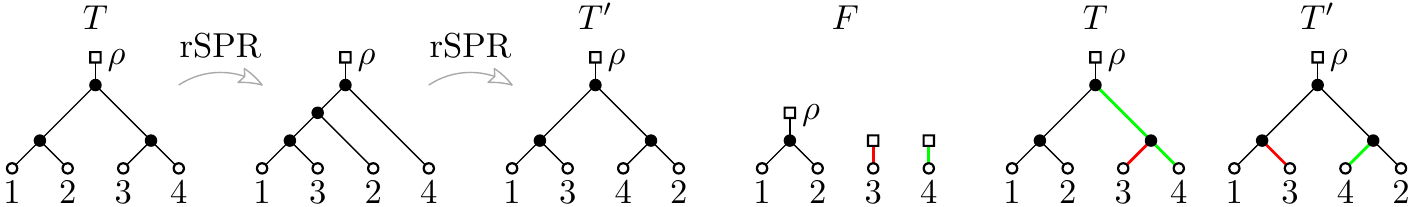}
  \caption{An rSPR-sequence of length two that transforms $T$ into $T'$, an agreement forest $F$ for
  $T$ and $T'$, and on the right embeddings of $F$ into $T$ and $T'$.}
  \label{fig:SNPR:MAF:treeTree}
\end{figure}

We start the section by informally describing the main ideas.
Consider the rSPR-sequence from $T$ to $T'$ shown in \cref{fig:SNPR:MAF:treeTree}.
This sequence first prunes and then regrafts the incoming edge of leaf 3, and then the incoming edge of leaf 4.
If we now look at this sequence and prune these edges again, but do not regraft them, then we obtain the forest $F$ shown in \cref{fig:SNPR:MAF:treeTree}.
The forest $F$ now represents the subtrees on which both $T$ and $T'$ ``agree''. 
Such $F$ is called an agreement forest for $T$ and $T'$ (defined precisely below).
In reverse and as also shown in \cref{fig:SNPR:MAF:treeTree}, $F$ can be embedded back into $T$ and $T'$ such that it covers all edges and vertices. 
The strength of such an agreement forest lies in the fact that it characterises the rSPR-distance of $T$ and $T'$ 
if it is optimal in some sense.

To generalise the idea of agreement forests to a tree and a network, we allow components that consist of a single edge. 
Intuitively, these components represent \SNPRP operations. 
We next make this precise and show how agreement forests can be used to characterise the SNPR-distance of $T$ and $N$. 

\begin{figure}[htb]
 \centering
 \includegraphics{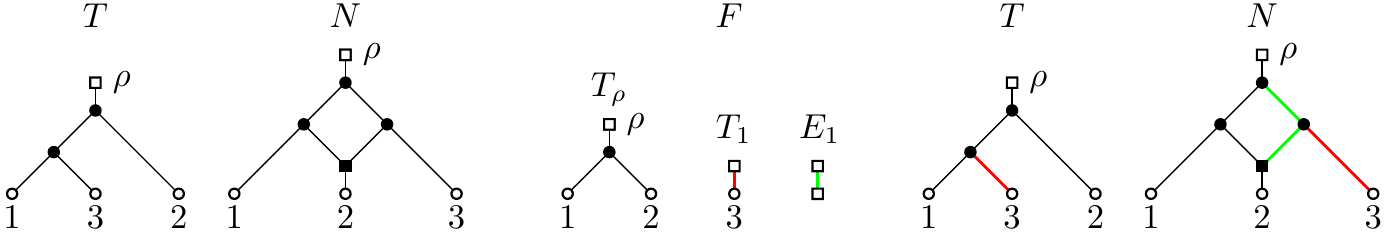}
 \caption{An agreement forest $F$ for the phylogenetic tree $T$ and the phylogenetic network
 $N$. On the right, embeddings of $F$ into $T$ and $N$.}
 \label{fig:SNPR:MAF:treeNetwork}   
\end{figure}

Let $T \in \trees$ and let $N \in \nets$ with $r$ reticulations.
For the purpose of the upcoming definition and much of this section, we regard the root $\root$ of
$T$ and $N$ as an element of the label sets $L(T)$ and $L(N)$, respectively. 
An \emph{agreement forest} $F$ for $T$ and $N$ is a collection
	$\{T_\root, T_1, T_2, \ldots, T_k, E_1, E_2, \ldots, E_r\}$, 
where $T_\root$ is an isolated vertex labelled $\root$, or a phylogenetic tree whose label set
includes $\root$, each $T_i$ with $i \in \{1, 2, \ldots, k\}$ is a phylogenetic tree,
and each $E_j$ with $j \in \{1, 2, \ldots, r\}$ is the graph that consists of a single directed edge
such that the following properties hold:
\begin{enumerate}[label=(\roman*)]
  \item The label sets $L(T_\root), L(T_1), L(T_2), \ldots, L(T_k)$ partition $\taxa \cup \{\root\}$.
  \item There exist simultaneous edge-disjoint embeddings of the trees 
   $$\{T_\root, T_1, T_2, \ldots, T_k\}$$ into $T$ that cover all edges of $T$.
  \item There exist simultaneous edge-disjoint embeddings of the graphs 
  	$$\{T_\root, T_1, T_2, \ldots, T_k, E_1, E_2, \ldots, E_r\}$$ into $N$ that cover all edges of $N$.
\end{enumerate}
Recall that ``cover'' here means that to each edge of $N$ an edge of a subdivision is mapped.
We refer to each element in $\{E_1, E_2, \ldots, E_r\}$ as a \emph{disagreement edge}.
To illustrate, \cref{fig:SNPR:MAF:treeNetwork} shows an agreement forest $F$ of a phylogenetic tree and a phylogenetic network.
We will show with \cref{clm:SNPR:MAF:treeNetwork} that an agreement forest for $T$ and $N$ always exists.

Let $F = \{T_\rho, T_1, T_2, \ldots, T_k, E_1, E_2, \ldots, E_r\}$ be an agreement forest for $T$ and $N$. 
Then $F$ is called a \emph{maximum agreement forest} for $T$ and $N$ if the number of
elements in the subset $\{T_\rho, T_1, T_2, \ldots, T_k\}$ of $F$ or, equivalently, $k$ is minimised.
Moreover, for $k$ being minimum, we set $m(T, N) = k + r = \abs{F} - 1$. 

Referring back to \cref{fig:SNPR:MAF:treeNetwork}, the agreement forest $F$ is a maximum agreement forest for $T$ and $N$.

For readers familiar with the notion of agreement forests for two phylogenetic trees $T$ and $T'$,
we note that the aforementioned definition of a maximum agreement forest coincides with its namesake concept for $T$ and $T'$ as introduced by Bordewich and Semple~\cite{BS05}.
The importance of the notion of maximum agreement forests for two phylogenetic trees lies in the following theorem.

\begin{theorem}[{{Bordewich and Semple \cite[Theorem 2.1]{BS05}}}]\label{clm:SNPR:MAF:trees}
Let $T, T' \in \trees$. Then $$\drSPR(T, T') = m(T, T')\text{.}$$
\end{theorem}

Next we show how the more general definition of agreement forests that is introduced in this paper
can be employed to characterise the SNPR-distance between a phylogenetic tree $T$ and a phylogenetic network $N$.  
We start with a `warm-up' for when $T$ is displayed by $N$.

\begin{lemma}\label{clm:SNPR:MAF:diplayedTree}
Let $N \in \nets$ with $r$ reticulations. Let $T \in D(N)$.
Then $$\dSNPR(T, N) = m(T, N) = r\text{.}$$
  \begin{proof}
    By \cref{clm:SNPR:shortestPath:displayedTree}, we have $\dSNPR(T, N) = r$ 
    and know that there exists an \SNPRP-sequence $\sigma = (T = N_0, N_1, \ldots, N_r = N)$ that transforms $T$ into $N$.
    Using $\sigma$, we now prove that $F = \{T = T_\rho, E_1, \ldots, E_r\}$ is an agreement forest for $T$ and $N$.
    The proof is by induction on $r$. 
    If $r = 0$, then $T = N$ and the claim trivially holds.
    Next, let $e$ be the edge added from $N_{i-1}$ to $N_i$ for $i = \{1, \ldots, r\}$.
    Note that $F_{i-1} = \{T, E_1, \ldots, E_{i-1}\}$ has an embedding into $N_i$ (as required for an agreement forest)
    that covers all edges except $e$.  
    Extending this embedding by mapping $E_i$ of $F_{i} = \{T, E_1, \ldots, E_{i}\}$ to $e$, 
    we get that $F_i$ is an agreement forest of $T$ and $N_i$.
    This is illustrated in \cref{fig:SNPR:MAF:treeToNetwork}.
    Hence, $F$ is an agreement forest for $T$ and $N$ and therefore
	$$r = \dSNPR(T, N) = \abs{F} - 1 \geq m(T, N)\text{.}$$ 
    
    \begin{figure}[htb]
	 \centering
	 \includegraphics{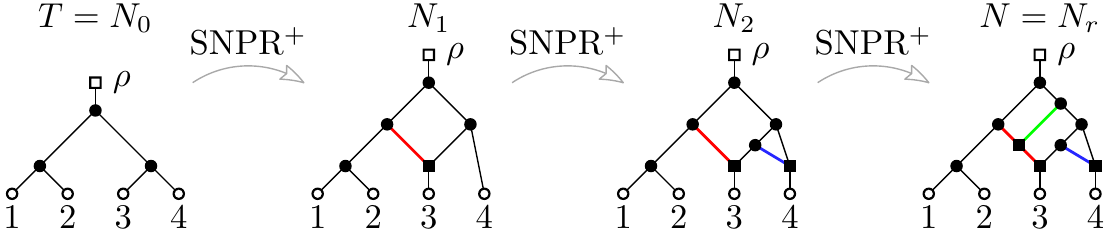}
	 \caption{An example of how to obtain an embedding into $N$ of an agreement forest 
	 $F = \{T, E_1, \ldots, E_r\}$ for $T$ and $N$ for the proof of \cref{clm:SNPR:MAF:diplayedTree}.}
	 \label{fig:SNPR:MAF:treeToNetwork}   
	\end{figure}
    
    To establish the other direction, let $F$ be a maximum agreement forest for $N$ and $T$.
    Recall that, by definition, $F$ contains $r$ disagreement edges and at least one other element.
    Thus, $$m(T, N) = \abs{F} - 1 \geq r + 1 - 1 = r = \dSNPR(T,N)\text{.}$$ 
    This completes the proof of the lemma.
  \end{proof}
\end{lemma}

We are now in a position to establish the main result of this section.

\begin{theorem}\label{clm:SNPR:MAF:treeNetwork}
Let $T \in \trees$, $N \in \nets$.
Then $$\dSNPR(T, N) = m(T, N)\text{.}$$
  \begin{proof}
    Let $r$ be the number of reticulations in $N$. 
    We first show that $m(T, N) \leq \dSNPR(T, N)$. 
    By \cref{clm:SNPR:shortestPath:nondisplayedTree}, there exists a phylogenetic tree $T'$
    that is displayed by $N$ such that 
    	$$\dSNPR(T, N) = \dSNPR(T, T') + \dSNPR(T', N) = \dSNPR(T, T') + r\text{.}$$  
    Hence, we have $m(T, T') = \dSNPR(T, T') = \dSNPR(T, N) - r$, where the first equality follows
    from \cref{clm:SNPR:isometric:trees} and \cref{clm:SNPR:MAF:trees}. 
    Moreover, by \cref{clm:SNPR:MAF:diplayedTree}, we have $m(T', N) = \dSNPR(T', N) = r$. 
    Let $F'$ be a maximum agreement forest for $T$ and $T'$, and let $F''$ be a maximum agreement
    forest for $T'$ and $N$. We know by \cref{clm:SNPR:MAF:diplayedTree} that such an $F''$ exists and that $T' \in F''$.
    Now, let 
    	$$F = F' \cup(F'' - \{T'\})\text{.}$$
	Since $F'$ has an embedding into $T'$, and $T'$ has an embedding into $N$,
	we get an embedding of $F'$ into $N$.
	This embedding covers all edges of $N$, except those to which the disagreement edges of $F''$ get mapped.
	Since $F$ contains both $F'$ and the disagreement edges of $F''$, it follows that $F$ is an agreement forest for $T$ and $N$. 
    Hence, $$m(T, N) \leq \abs{F} - 1 = \abs{F'} + \abs{F''} - 2 = \dSNPR(T, T') + \dSNPR(T',N) = \dSNPR(T,N)\text{.}$$

	We next show that $\dSNPR(T, N) \leq m(T, N)$. 
	The proof is by induction on the size $\abs{F}$ of a maximum agreement forest $F$ for $T$ and $N$,
	which we can write as $$F = \{T_\rho, T_1, T_2,\ldots, T_k, E_1,, E_2,\ldots, E_r\}\text{.}$$
	If $\abs{F} = 1$, that is $F = \{T\}$, then $N = T$ and, so, $\dSNPR(T, N) = 0$. 
	Now assume that the inequality holds for all pairs of a phylogenetic tree and a phylogenetic
	network on the same leaf set for which there exists a maximum agreement forest whose number of components is at most $k + r$. 
	If $r = 0$, then $N$ is a phylogenetic tree and $F = \{T_\rho, T_1, T_2, \ldots, T_k\}$.
	Then it follows from \cref{clm:SNPR:MAF:trees} that $\drSPR(T, N) \leq m(T, N)$. 
	Moreover, by \cref{clm:SNPR:isometric:trees}, we have that $\dSNPR(T, N) = \drSPR(T, N) \leq m(T, N)$.
	
	We may therefore assume that $r > 0$. 
	Let $v$ be a reticulation in $N$ that has no reticulation as an ancestor. 
	For each component $C_i \in F$, let $\epsilon(C_i)$ be the set of edges in $N$ that is used to embed $C_i$ into $N$ such that 
	$$\cE = \{\epsilon(T_\rho), \epsilon(T_1), \epsilon(T_2), \ldots, \epsilon(T_k), \epsilon(E_1),	\epsilon(E_2), \ldots, \epsilon(E_r)\}$$ 
	is a partition of the edge set of $N$. 
	Since $F$ is an agreement forest for $T$ and $N$, such a partition exists.
	
	Now, let $(u, v)$ and $(u', v)$ be the reticulation edges incident with $v$. 
	Without loss of	generality, we may assume that $(u, v) \in \epsilon(E_i)$ for some $i \in \{1, 2, \ldots, r\}$.
	Note that if $(u', v) \in \epsilon(T_j)$ for some $j \in \{\rho, 1, 2, \ldots, k\}$ (i.e., no disagreement edge is mapped to $(u', v)$),
		then $\epsilon(T_j)$ also contains the outgoing edge of $v$.
	Otherwise, if $(u', v) \in \epsilon(E_j)$ for some $j \in \{1, 2, \ldots, r\}$, $j \neq i$, 
	then we may assume without loss of generality that $\epsilon(E_j)$ (and not $\epsilon(E_i)$) contains the outgoing edge of $v$. 
	Let $F' = F - \{E_i\}$, and let $N'$ be the phylogenetic network obtained from $N$ be deleting $(u, v)$ 
	and suppressing the resulting two degree-2 vertices. 
	We next show that $F'$ is an agreement forest for $T$ and $N'$. 
	By the choice of $v$, recall that $u$ is a tree vertex. 
	Let $w$ be the second child of $u$ and let $C_j$ be the component in $F$ such that $(u, w) \in \epsilon(C_j)$. 
	(Note that if $N$ contains a parallel edge such that $u = u'$ then $w = v$.)
	Set $\epsilon'(C_j) = \epsilon(C_j) \cup (\epsilon(E_i) - \{(u, v)\})$. 
	This is illustrated in \cref{fig:SNPR:MAF:embeddingChange}. 
	Note that $\epsilon'(C_j) = \epsilon(C_j)$ precisely if $\epsilon(E_i) = \{(u, v)\}$. 
	As $\cE$ is	an embedding of $F$ into $N$ that partitions the edge set of $N$,
		$$\cE' = (\cE - \{\epsilon(C_j), \epsilon(E_i)\}) \cup \{\epsilon'(C_j)\}$$ 
	partitions the edge set	of $N'$ and induces an embedding of $F'$ in $N'$. 
	Hence, $F'$ is an agreement forest for $N'$ and $T'$. 
	Since $\abs{F'} < \abs{F}$, it now follows from the induction hypothesis that there
	exists an SNPR-sequence from $T$ to $N'$ whose length is at most $\abs{F'} = k + r - 1$. 
	Furthermore, by	construction, $N$ can be obtained from $N'$ by a single \SNPRP.
	Taken together, this implies that 
	$$\dSNPR(T, N) \leq \dSNPR(T, N') + 1 \leq k + r - 1 + 1 = m(T,N)\text{.}$$
	
	\begin{figure}[htb]
	  \centering
	  \includegraphics{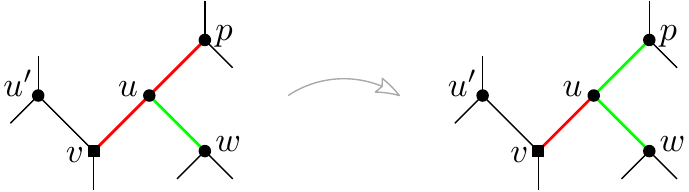}
	  \caption{An example for the proof of \cref{clm:SNPR:MAF:treeNetwork} showing how $E_i$ (red) and $C_j$ (green) embed into $N$. 
	  On the left, $\epsilon(E_i) = \{(p, u), (u, v)\}$ and $\epsilon(C_j)$ contains (at least) $(u, w)$. 
	  Thus, on the right, $\epsilon'(C_j)$ is obtained from $\epsilon(C_j)$ by adding $(p, u)$.}
	  \label{fig:SNPR:MAF:embeddingChange}
	\end{figure}
	
	Combining both inequalities establishes the theorem.
  \end{proof}
\end{theorem}

\section{Properties of shortest \SNPR-sequences connecting two networks}
\label{sec:shortestPaths}

In this section, we analyse properties of shortest SNPR-sequences that connect a pair of
phylogenetic networks and investigate whether or not the three classes of tree-child,
reticulation-visible, and tree-based networks isometrically embed into the class of all phylogenetic networks.
We start with some definitions that are used throughout this section.
For any non-negative integer $r$, \emph{tier} $r$ of $\nets$ is the subset of networks in $\nets$
that have exactly $r$ reticulations. Note that tier $0$ equals $\trees$.
For $N, N' \in \nets$, let $\sigma = (N = N_0, N_1, \ldots, N_k = N')$ be an SNPR-sequence from $N$ to $N'$.
We say that $\sigma$ \emph{horizontally traverses} tier $r$ if $\sigma$ contains two networks
$N_{i-1}$ and $N_{i}$ with $i \in \{1, 2, \ldots, k\}$ such that both have $r$ reticulations; 
i.e., $N_{i}$ can be obtained from $N_{i-1}$ by a single \SNPRZ.

Let $N, N' \in \nets$ with $r$ and $r'$ reticulations, respectively. 
Without loss of generality, we may assume that $r \leq r'$.
From a computational viewpoint and in trying to shrink the search space when computing $\dSNPR(N, N')$, 
	it would be desirable if there always exists a shortest SNPR-sequence connecting $N$ and $N'$ that traverses exactly one tier horizontally.
In particular, if $r < r'$ it would have positive implications for computing $\dSNPR(N, N')$ if all
\SNPRZ operations could be pushed to be the beginning or the end of a shortest SNPR-sequence for $N$ and $N'$.
On the other hand, if $r = r'$, then the existence of a shortest SNPR-sequence from $N$ to $N'$
whose networks all belong to tier $r$ would allow us to compute $\dSNPR(N, N')$ by considering only tier $r$.
In what follows, we present several results showing that the existence of a shortest SNPR-sequence
with such properties cannot be guaranteed. 
For each result we provide a small counterexample or a family of counterexamples. 
Furthermore, the networks in these examples can be extended to contain more reticulations and taxa.
See also the discussion at the end of this section.

%
%
\begin{lemma}\label{clm:SNPR:shortestPath:zpz}
Let $n \geq 4$.
Let $N, N' \in \nets$ with $r$ and $r'$ reticulations, respectively, such that $r < r'$.
Then there does not necessarily exist a shortest \SNPR-sequence from $N$ to $N'$ that
traverses at most one tier horizontally.
  \begin{proof}
    To prove the statement, we show that every shortest SNPR-sequence for the two
    phylogenetic networks $N$ and $N'$ that are depicted in \cref{fig:SNPR:shortestPath:zpz}
    traverses at least two tiers horizontally.  
    
	We start by observing four differences between $N$ and $N'$:
	\begin{enumerate}[label=(\arabic*)]
  	  \item Leaf 1 is a descendant of a reticulation in $N$, but not in $N'$. 
  	  \item Leaves 1 and 4 form a cherry in $N'$, but not in $N$.  
  	  \item Leaves 2 and 3 form a cherry in $N'$, but not in $N$.
  	  \item Leaves 2 and 3 are descendants of two reticulations in $N'$, but not in $N$.  	  
	\end{enumerate}

    \begin{figure}[htb]
	 \centering
	 \includegraphics{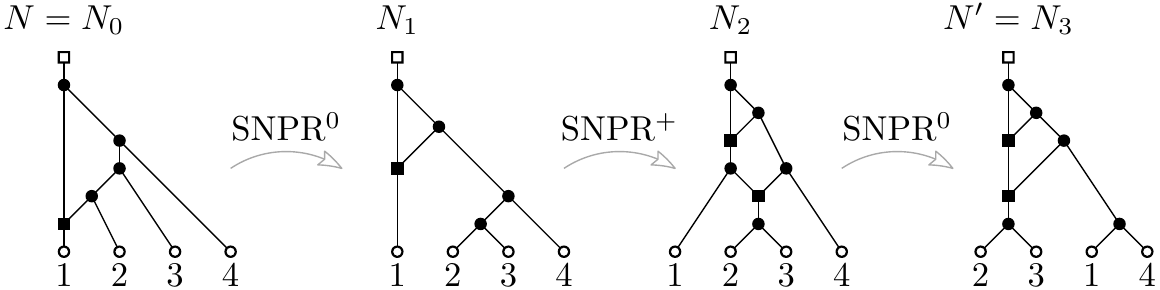}
	 \caption{For the two networks $N$ and $N'$ shown, every shortest \SNPR-sequence between them traverses two tiers horizontally.}
	 \label{fig:SNPR:shortestPath:zpz}
	\end{figure}

	Since $N'$ has one more reticulation than $N$, at least one \SNPRP is required to	transform $N$ into $N'$.
	Also note that an \SNPRP cannot in general create a cherry.
	Furthermore, note that an \SNPRZ on $N$ (or a network derived from $N$ by an \SNPRP) can create at most one cherry. 
	Therefore, to transform $N$ into $N'$ at least three \SNPR operations are necessary and thus $\dSNPR(N, N') > 2$.
	Consequently, referring back to the networks shown in \cref{fig:SNPR:shortestPath:zpz},
	$$\sigma = (N = N_0, N_1, N_2, N_3 = N')$$ is a shortest SNPR-sequence from $N$ to $N'$ that
	horizontally traverses tier 1 and tier 2.
	
	To establish the statement, it is therefore sufficient to show that there exists no SNPR-sequence,
	say $$\sigma^* = (N, M, M', N')\text{,}$$ such that $M$ can be obtained from $N$ by an \SNPRP, or $N'$ can be obtained from $M'$ by an \SNPRP.
	Note that a sequence that uses an \SNPRP (or an \SNPRM) to transform $M$ into $M'$ would either be covered by one of these two cases
	or would be a sequence that traverses two tiers horizontally like $\sigma$.
	We thus proceed by distinguishing the first two cases.
	
	First, assume that $\sigma^*$ exists and that $M$ has been obtained from $N$ by an \SNPRP. 
	Then $M$ and $N'$ have the same four differences as listed above for $N$ and $N'$ with the
	exception that either leaf 2 or 3 (but not both) is possibly a descendant of two reticulations in
	$M$. Suppose that $M$ is indeed obtained from $N$ by 
		(i) subdividing the edge directed into 1 with a new vertex $u$, subdividing the edge directed into
		2 with a new vertex $v$, and adding the new edge $(u,v)$, 
	or 
		(ii) subdividing the edge directed into 1 with a new vertex $u$, subdividing the edge directed
		into 3 with a new vertex $v$, and adding the new edge $(u,v)$. 
	Then $M$ would equal either the network $M_1$ or $M_2$ shown in \cref{fig:SNPR:shortestPath:zpzMs}.
	In both cases, it requires two SNPR operations to transform $M$ into
	a network, say $M^*$, 
	in which leaf 1 is not a descendant of any reticulation and leaves 2 and 3 are descendants of two reticulations.
	One such $M^*$ is shown in \cref{fig:SNPR:shortestPath:zpzMs}. 
	However, $M^* \neq N'$ and, so, it would take in total at least three SNPR operations to transform $M$ into $N'$. 
	Now, suppose that $M$ is obtained from $N$ by an $\SNPR^+$ other than (i) or (ii). 
	With similar observations as above we note that again at least three SNPR operations are necessary to transform $M$ into $N'$.
	Hence, we conclude that $M$ has not been obtained from $N$ by an \SNPRP.

    \begin{figure}[htb]
	 \centering
	 \includegraphics{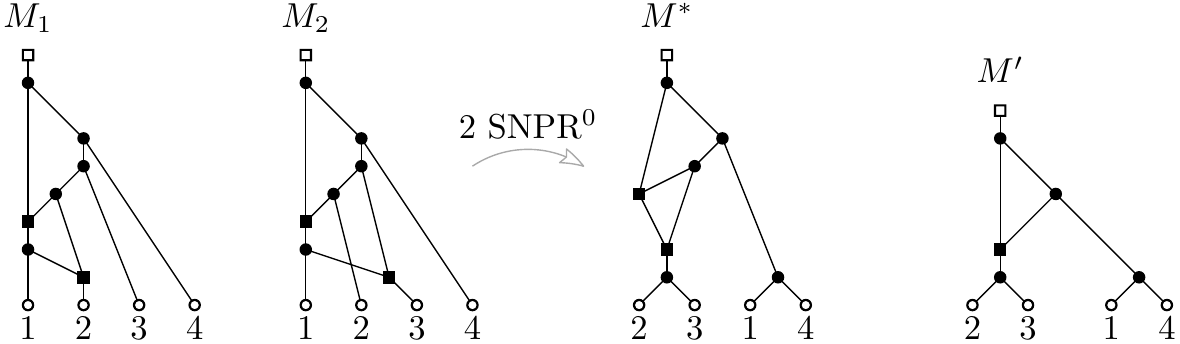}
	 \caption{Networks in SNPR-sequences from $N$ to $N'$ of \cref{fig:SNPR:shortestPath:zpz} for the proof of \cref{clm:SNPR:shortestPath:zpz}.}
	 \label{fig:SNPR:shortestPath:zpzMs}
	\end{figure}

	Second, assume that $\sigma^*$ exists and that $N'$ has been obtained from $M'$ by an \SNPRP or,
	equivalently, $M'$ has been obtained from $N'$ by an \SNPRM. 
	Then $M'$ is as shown in \cref{fig:SNPR:shortestPath:zpzMs} since each of the three \SNPRM operations that can be
	applied to $N'$ results in the same network $M'$. 
	Because of the aforementioned differences between $N$ and $N'$ that are also differences
	between $N$ and $M'$ with the exception that 2 and 3 are descendants of only a single
	reticulation in $M'$, it takes at least three SNPR operations to transform $N$ into $M'$.
	Consequently, $N'$ has not been obtained from $M'$ in $\sigma^*$ by an \SNPRP.
	
	Lastly, since neither $M$ nor $N'$ has been obtained from $N$ and $M'$, respectively, by an
	\SNPRP, it follows that $\sigma^*$ cannot be chosen so that no tier is horizontally traversed.
	This completes the proof. 
  \end{proof}
\end{lemma}

We next shows that, for two phylogenetic networks $N$ and $N'$ that both have $r$ reticulations,
every shortest SNPR-sequence from $N$ to $N'$ may contain a phylogenetic tree. 
Hence, to compute $\dSNPR(N, N')$ it may be necessary to search in the space of all phylogenetic networks with at most $r$ reticulations.

\begin{lemma}\label{clm:SNPR:shortestPath:down}
Let $r \geq 2$ and $n \geq 2r + 2$.
There exist $\bar N_r, \bar N_r' \in \nets$ with $r$ reticulations such that every shortest
\SNPR-sequence from $\bar N_r$ to $\bar N_r'$ contains a phylogenetic tree.
  \begin{proof}
    To prove the statement, we show that every shortest SNPR-sequence 
		$$\sigma = (\bar N_r = N_0, N_1, \ldots, N_k = \bar N_r')$$
	connecting the two phylogenetic networks $\bar N_r$ and $\bar N_r'$ depicted in
	\cref{fig:SNPR:shortestPath:down} has length $2k$, 
	for each $i \in \{1, 2, \ldots, r\}$, $N_i$ is obtained from $N_{i - 1}$ by an \SNPRM and for each
	$i \in \{r + 1, r + 2, \ldots, 2r\}$, $N_i$ is obtained from $N_{i - 1}$ by an \SNPRP. 
	Since $\bar N_r$ and $\bar N_r'$ both have $r$ reticulations, this implies that $\sigma$ contains
	a phylogenetic tree. 
	Note that $\sigma$ exists because we can transform $\bar N_r$ into $\bar N_r'$ by removing each
	reticulation edge in $\{e_1, e_2, \ldots, e_r\}$ with an \SNPRM and then adding each edge $\{e_1',
	e_2', \ldots, e_r'\}$ with an \SNPRP.

    \begin{figure}[htb]
	 \centering
	 \includegraphics{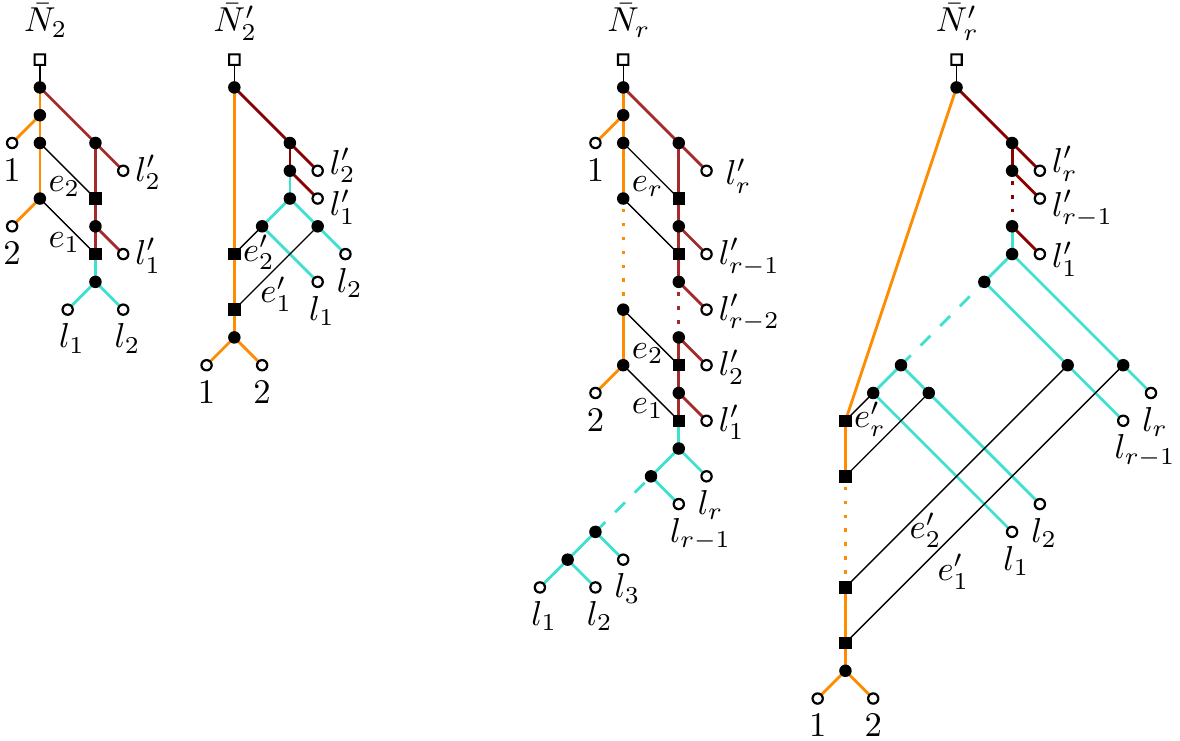}
	 \caption{Construction that is used in the proof of \cref{clm:SNPR:shortestPath:down} to show
	 that, for each $r \geq 2$, there exist two phylogenetic networks $\bar N_r$ and $\bar N_r'$ such
	 that every shortest \SNPR-sequence from $\bar N_r$ to $\bar N_r'$ contains a
	 phylogenetic tree.} 
	 \label{fig:SNPR:shortestPath:down}
	\end{figure}

	We pause to observe three properties of $\bar N_r'$ that will be crucial for the remainder of
	this proof: 
	\begin{enumerate}[label=(P\arabic*)]
	  \item	For each $i \in \{1, 2, \ldots, r\}$, the leaf $l_i$ is a sibling of a
	  reticulation.
	  \item Leaves 1 and 2 form a cherry, and descendants of all reticulations.
	  \item	There exists a directed path $(\root, w, v_1, v_2, \ldots, v_r)$, 
	  where $\root$ is the root, $w$ is the child of $\rho$, and each $v_i$ with $i \in \{1, 2, \ldots, r\}$ is a
	  reticulation.
	\end{enumerate}
	To illustrate, for $r = 2$, the networks $\bar N_2$ and $\bar N_2'$ are shown in
	\cref{fig:SNPR:shortestPath:down}.
	
	Now assume that there exists an SNPR-sequence
		$$\sigma^* = (\bar N_r = M_0, M_1, M_2, \ldots,M_{k'} = \bar N_r')$$ 
	from $\bar N_r$ to $\bar N_r'$ of length $k' \leq 2r$ that is distinct from $\sigma$.
	Let 
		$$O^* = (o_1, o_2, \ldots, o_{k'})$$ 
	be the sequence obtained from $\sigma^*$ such that for each	$i \in \{1, 2, \ldots, k'\}$ the following holds:
	\begin{itemize}
	\item $o_i = 0$ if $M_i$ is obtained from $M_{i-1}$ by an \SNPRZ,
	\item $o_i = +$ if $M_i$ is obtained from $M_{i-1}$ by an \SNPRP, or
	\item $o_i = -$ if $M_i$ is obtained from $M_{i-1}$ by an \SNPRM.
	\end{itemize}
	Let $m$ be the number of elements in $O^*$ that are equal to $-$.
	
	\noindent{\textbf{Case 1.}} 
	Assume that $m > r$. 
	Since $\bar N_r$ and $\bar N_r'$ both have $r$ reticulations, $O^*$ contains exactly $m$ elements
	that are equal to $+$. Hence, $k' \geq 2m > 2r$; a contradiction.
	
	\noindent{\textbf{Case 2.}} 
	Assume that $m < r$.
	Again, since $\bar N_r$ and $\bar N_r'$ both have $r$ reticulations, $O^*$ contains exactly $m$
	elements that are equal to $+$.
	Thus, with $k' \leq 2r$, it follows that $O^*$ contains at most $2(r - m)$ elements that are equal
	to $0$.
	Let $i$ be an element in $\{1, 2, \ldots, k'\}$ such that $o_i = +$. 
	Then, the number of leaves in $\{l_1, l_2, \ldots, l_r\}$ that are siblings of different reticulations in
	$M_{i-1}$ and $M_i$ differs by at most one. 
	Therefore, we need at least $k_1 \geq r - m$ \SNPRZ operations to obtain a network from $\bar N_r$
	that satisfies (P1). 
	Similarly, the number of vertices on a directed path that consists only of reticulations in
	$M_{i-1}$ and $M_i$ differs by at most one. 
	Therefore, we need at least $k_2 \geq r - m$ \SNPRZ operations to obtain a network from $\bar N_r$
	that satisfies (P3).
	
	Let $i \in \{1, 2, \ldots, k'\}$ such that $o_i = 0$.
	Assume that the number of leaves in $\{l_1, l_2,$ $\ldots, l_r\}$ that are siblings of reticulations in
	$M_i$ is greater than this number in $M_{i-1}$.
	Then, the \SNPRZ operation to obtain $M_i$ from $M_i$ either regrafts such a leaf $l_j$ as sibling to the incoming edge of a reticulation
	 or regrafts a reticulation edge to the incoming edge of such a leaf. 
	Therefore this operation cannot have increased the number of vertices that lie on a directed path 
	of reticulations in $M_i$ compared to $M_{i-1}$.
	Similarly, if the number of vertices that lie on a directed path of reticulations in $M_i$ is
	greater than that number in $M_{i-1}$, then the number of leaves in $\{l_1, l_2, \ldots, l_r\}$
	that are siblings of reticulations is not greater in $M_i$ than in $M_{i-1}$.
	Again, an \SNPRZ operation cannot change both values for these networks at the same time.
	Overall, we observe that the $k_1$ \SNPRZ used to satisfy property (P1) affect the leaves $l_j$ and
	reticulation edges, whereas the $k_2$ \SNPRZ used to satisfy property (P3) affect the leaves
	$l_j'$ and (possibly) leaf 1. 
	It follows that $k_1 = k_2 = (r - m)$ and, so, $k' = 2r$.
	
	Lastly, to see that $M_{k'}$ does not satisfy property (P2), observe that neither the $k_1 + k_2$
	\SNPRZ operations nor the $2m$ \SNPRM and \SNPRP operations that are used to satisfy
	(P1) and (P3) result in a network that simultaneously satisfies (P2).
	Hence, it follows that at least one additional \SNPRZ is needed to transform $\bar N_r$ into $\bar
	N_r'$; thereby contradicting that $k' \leq 2r$.
	
	\noindent{\textbf{Case 3.}} Assume that $m = r$. 
	Since $\bar N_r$ and $\bar N_r'$ both have $r$ reticulations and $k'\leq 2r$, it follows that $k' = 2r$.
	We complete the proof by showing that, for each $i \in \{1, 2, \ldots, r\}$, we have $o_i = -$ and,
	for each $i \in \{r + 1, r + 2, \ldots, 2r\}$, we have $o_i = +$.
	Assume that, for some $i \leq r$, we have $o_i = +$. 
	Choose $i$ to be as small as possible.
	Let $v$ be the unique reticulation in $M_i$ that is not a reticulation in $M_{i-1}$. 
	Then $v$ does not have leaves 1 and 2 as descendants and a leaf in $\{l_1, l_2, \ldots, l_r\}$ as a
	sibling of a reticulation. 
	Now, as $O^*$ does not contain an element equal to 0, there exists an element $o_j = -$ with $j >
	i$ such that $M_j$ does not contain the reticulation edge that was added in transforming 
	$M_{i-1}$ into $M_i$. 
	In turn, this implies that the remaining $r-1$ \SNPRP cannot transform $\bar N_r$ into a network
	that satisfies (P1) and (P3). Hence, if $m = r$, then $\sigma^*=\sigma$.
	
	Combining all three cases establishes the statement.
  \end{proof}
\end{lemma}

\begin{figure}[htb]
 \centering
 \includegraphics{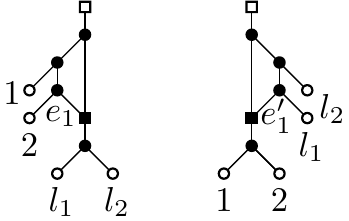}
 \caption{Two phylogenetic networks with one reticulation such that every shortest \SNPR-sequence
 connecting them contains a phylogenetic tree.}
 \label{fig:SNPR:shortestPath:downOne}
\end{figure}

Recall that the statement of \cref{clm:SNPR:shortestPath:down} requires $\bar N_r$ and
$\bar N_r'$ to have at least two reticulations. 
Using a slightly different construction than that for $\bar N_r$ and $\bar N_r'$, 
\cref{fig:SNPR:shortestPath:downOne} shows two phylogenetic networks that both have one reticulation
such that every shortest SNPR-sequence connecting these two networks contains a phylogenetic tree. 
While omitting a formal proof, we note that it can be checked by following the same ideas as in the proof of \cref{clm:SNPR:shortestPath:down}.

Bordewich et al.~\cite[Proposition 7.5]{BLS17} showed that  
$$\dSNPR(N, N') \leq \min\{\dSNPR(T, T') \colon T \in D(N) \textnormal{ and } T' \in D(N')\} + r +
r'\text{,}$$ where $N, N' \in \nets$ with $r$ and $r'$ reticulations, respectively. 
\cref{clm:SNPR:shortestPath:down} implies that this upper bound is sharp, 
for example, for the networks $\bar N_r$ and $\bar N_r'$ in \cref{fig:SNPR:shortestPath:down}. 

The next lemma shows that, for two phylogenetic networks $N$ and $N'$ that both have $r$ reticulations, 
every shortest SNPR-sequence from $N$ to $N'$ may contain a network that has more than $r$ reticulations. 
In particular, to compute $\dSNPR(N, N')$ it is not sufficient to only search the space of all phylogenetic networks that have at most $r$ reticulations.

\begin{lemma}\label{clm:SNPR:shortestPath:up}
Let $n \geq 2$, $r \geq 3$, and let $N, N' \in \nets$ with $r$ reticulations.\\
There does not necessarily exist a shortest \SNPR-sequence from $N$ to $N'$ such
that each network in the sequence has at most $r$ reticulations.
  \begin{proof}
    To establish the lemma, we show that every shortest SNPR-sequence that connects the two phylogenetic networks $N$ and $N'$ 
    as depicted in \cref{fig:SNPR:shortestPath:up} contains a network with four reticulations. 
    First observe that $\dSNPR(N, N') \geq 2$ and, so, the SNPR-sequence $(N, N_1, N')$ is of minimum length.

    \begin{figure}[htb]
	 \centering
	 \includegraphics{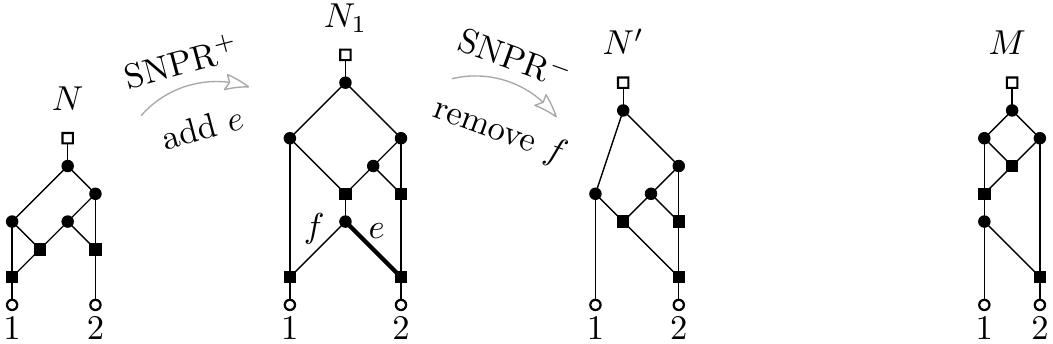} 
	 \caption{Example that is used in the proof of \cref{clm:SNPR:shortestPath:up}
	 	 showing two networks $N$ and $N'$, which have three reticulations each 
	 	 and for which every shortest \SNPR-sequence between them contains a network with four reticulations.}
	 \label{fig:SNPR:shortestPath:up} 
	\end{figure}

	We now show that there exists no SNPR-sequence $(N, M, N')$
	such that $M$ is obtained from $N$ by an \SNPRM or \SNPRZ.
	Towards a contradiction, assume that $M$ is obtained from $N$ by an \SNPRM. 
	Clearly, leaf 1 is a child of a reticulation in $M$.
	Moreover, as $M$ has two reticulations,	it follows that $N'$ is obtained from $M$ by an \SNPRP
	and that leaf 1 is still a child of a reticulation in $N'$; a contradiction. 
	Now assume that $M$ is obtained from $N$ by an \SNPRZ. 
	If leaf 1 is a child of a reticulation in $M$, then $\dSNPR(M, N') > 1$. 
	We may therefore assume that leaf 1 is not a child of a	reticulation in $M$. 
	Hence, $M$ is the network that is shown on the right-hand side of \cref{fig:SNPR:shortestPath:up}.
	Note that in $M$ (contrary to $N'$) the leaf 1 is descendant of a reticulation and all three reticulations are on a directed path.
	We observe that changing either of these properties with a single \SNPRZ cannot change the other property.
	Therefore $\dSNPR(M, N') > 1$; again a contradiction.
  \end{proof}
\end{lemma}

The next theorem is an immediate consequence of \cref{clm:SNPR:shortestPath:down,clm:SNPR:shortestPath:up} and \cref{fig:SNPR:shortestPath:downOne}.

\begin{theorem}\label{clm:SNPR:nonIsomorphic:tiers}
Let $\cC_r$ be the class of all phylogenetic networks in $\nets$ that have $r$ reticulations.
If $n \geq 4$ and $r \geq 1$, then $\cC_r$ does not isometrically embed into the class of all phylogenetic networks $\nets$. 
Moreover, if $n \geq 2$ and $r \geq 3$, then $\cC_r$ does not isometrically embed 
into the class of all phylogenetic networks in $\nets$ with at most $r$ reticulations.
\end{theorem}
 
We now consider different classes of phylogenetic networks and ask if they isometrically
embed into the class of all phylogenetic networks.  
As we will see, we answer this question negatively for tree-child networks $\tchinets$,
reticulation-visible networks $\retvisnets$, and tree-based networks $\tbasednets$.

\begin{figure}[htb]
	 \centering
	 \includegraphics{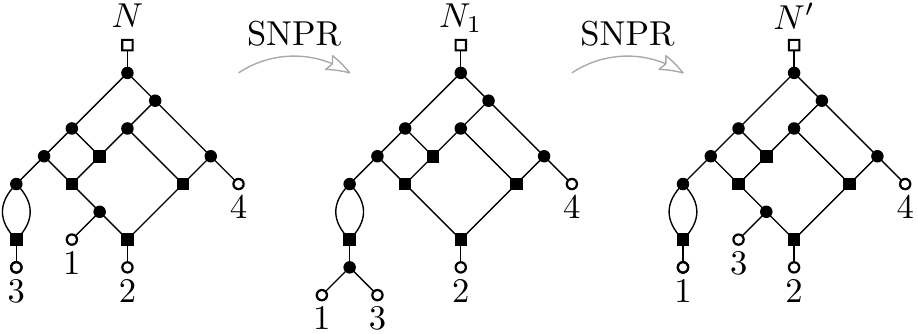}
	 \caption{
	 Example that is used in the proof of \cref{clm:SNPR:nonIsomorphic:classes} 
	 to show that $\tbasednets$ does not embed isometrically into $\nets$.
	 }
	 \label{fig:SNPR:nonIsomorphic:tbased}
	\end{figure}

\begin{proposition}\label{clm:SNPR:nonIsomorphic:classes}
Let $\class \in \{\tchinets, \retvisnets, \tbasednets\}$ with $n \geq 4$.\\
Then $\class$ does not embed isometrically into $\nets$ under \SNPR.
  \begin{proof}
    To establish the theorem, we give explicit examples of two networks $N$ and $N'$ that are in
    $\class$ such that $d_{\SNPR_{\class}}(N, N') > d_{\SNPR_{\nets}}(N, N')$.

	\begin{figure}[htb]
	 \centering
	 \includegraphics{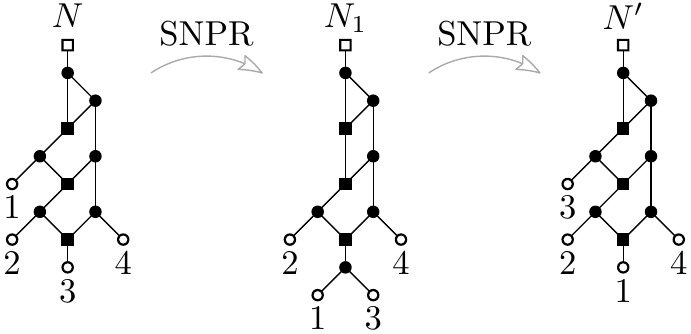}
	 \caption{Example that is used in the proof of \cref{clm:SNPR:nonIsomorphic:classes} 
	 to show that neither $\tchinets$ nor $\retvisnets$ embeds isometrically into
	 $\nets$. 
	}
	 \label{fig:SNPR:nonIsomorphic:tchi} 
	\end{figure}

 	Let $\class = \tchinets$ (resp. $\class = \tbasednets$). 
 	Consider the two tree-child (resp. tree-based) networks $N$ and $N'$ that are shown in
 	\cref{fig:SNPR:nonIsomorphic:tchi} (resp. \cref{fig:SNPR:nonIsomorphic:tbased}). 
 	Then $\sigma = (N, N_1, N')$ is an SNPR-sequence for $N$ and $N'$. 
 	Note that $N'$ can be obtained from $N$ by swapping the labels $1$ and $3$.
 	Since leaf $3$ is the child of a reticulation in $N$, it cannot be pruned with an \SNPRZ in $N$.
 	The sequence $\sigma$ thus prunes the edge incident to leaf $1$ to regraft it above leaf $3$, 
 	which then enables the edge incident to leaf $3$ to be pruned and regrafted to the former position of leaf $1$.
 	
 	Towards a contradiction, assume that there exists an SNPR-sequence $\sigma^* = (N, M, N')$ distinct from $\sigma$.
 	Suppose $\sigma^*$ does not start by pruning the edge incident to leaf $1$. 
 	Then leaf $1$ has to be moved from $M$ to $N'$.
 	Furthermore, the edge incident to leaf $3$ cannot be pruned in $N$, so leaf 3 has to be moved from $M$ to $N'$. 
 	However, making both these changes is not possible with a single SNPR operation.
 	Therefore, $\sigma$ is the unique SNPR-sequence in $\nets$ of length two that connects $N$ and $N'$. 
 	Hence, as $M$ is not tree child (resp. tree based), we have 
 		$$d_{\SNPR_{\class}}(N, N') > d_{\SNPR_{\nets}}(N, N') = 2\text{.}$$ 
 	Noting that $M$ in \cref{clm:SNPR:nonIsomorphic:classes} is not reticulation visible, the
 	same argument holds for when $\class = \retvisnets$.
  \end{proof}
\end{proposition}

While Francis and Steel~\cite{FS15} allow tree-based networks to have edges in parallel,
we can also show that the class of all tree-based networks without parallel edges on $n$ leaves is not isometrically embedded into the
class of all phylogenetic networks on $n$ leaves either. 
For this we reuse the proof of \cref{clm:SNPR:nonIsomorphic:classes} with a counterexample obtained 
by subdividing each of the two edges in parallel that are shown in \cref{fig:SNPR:nonIsomorphic:tbased} with a new vertex, 
say $u$ and $v$, and adding a new edge $(u, v)$.

Lastly, the networks presented in this section may seem rather small. 
However, they can be regarded as skeletons of larger networks with the same properties. 
For instance, in all examples that we used to establish the results of this section, leaves
can be replaced with subtrees and subnetworks. 
Furthermore, some edges can be subdivided to add further reticulation edges or subtrees to obtain
larger networks with the same properties.

\section{Concluding remarks}
\label{sec:discussion}

In this paper, we have established the first results related to calculating the SNPR-distance, which
is an NP-hard problem.
In the first part, we have considered the special case of computing this distance between a
phylogenetic tree $T$ and a phylogenetic network $N$. 
In this particular case, computing the SNPR-distance is fixed-parameter tractable when parameterised
by this distance and can be calculated by solving several instances of the rSPR-distance problem.
Additionally, we have characterised the SNPR-distance of $T$ and $N$ in terms of agreement forests.
This result lends itself to an algorithm that works directly on $T$ and $N$ without having to solve
multiple instances of the rSPR-distance problem between two trees.
In the second part, we have turned to the SNPR-distance problem between two
phylogenetic networks $N$ and $N'$ and presented several results on shortest SNPR-sequences for $N$
and $N'$ with $r$ and $r'$ reticulations, respectively. 
These results show that the search space for computing the SNPR-distance of $N$ and $N'$ cannot in
general be pruned to networks whose number of reticulations is at least $\min\{r, r'\}$ or at most
$\max\{r, r'\}$. Furthermore, if $N$ and $N'$ are both tree child, reticulation visible, or tree
based, the search space cannot in general be restricted to these network classes.

As alluded to in the introduction, Gambette et al.~\cite{GvIJLPS17} have introduced a slightly
different operation that generalises rSPR to phylogenetic networks. 
The main difference between their operation and SNPR is that they allow for an additional operation
which is called a head move.
In the language of this paper, let $N$ be a phylogenetic network, and let $(u, v)$ be an edge of
$N$ such that $v$ is a reticulation. 
Then, the operation of deleting $(u, v)$, suppressing $u$, subdividing an edge that is not an
ancestor of $v$ with a new vertex $u'$, and adding the edge $(v,u')$ is referred to as a \emph{head
move}. Interestingly, if we generalise the SNPR operation by, additionally, allowing for head moves,
the properties of shortest SNPR-sequences that we have revealed in \cref{sec:shortestPaths} and
that may appear to be undesirable with regards to practical search algorithm do not change.
On the positive side, a characterisation of the SNPR-distance between a phylogenetic tree and a
phylogenetic network in terms of agreement forest is possible and a result equivalent to
\cref{clm:SNPR:MAF:treeNetwork} can be established. 
For further details, we refer the interested reader to the first author's PhD
thesis~\cite{SpOPhyN} which establishes results equivalent to the ones presented in this paper for when
one allows for head moves.

\begin{figure}[htb]
\centering
  \includegraphics{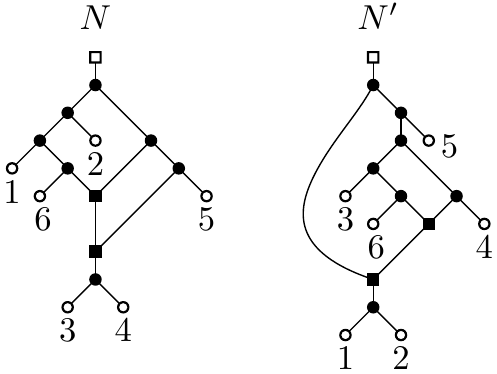}
  \caption{Two phylogenetic networks $N$ and $N'$ for which every shortest SNPR-sequence prunes at
  least one edge twice.}
  \label{fig:SNPR:doublePrune}
\end{figure}

We close this paper by asking whether the notion of agreement forests can be further generalised
to computing the SNPR-distance between two phylogenetic networks, regardless of whether head
moves are also allowed.
As mentioned above, a shortest sequence between $N$ and $N'$ might
have to traverse a tier with more or less reticulations than $N$ and $N'$. 
It is unclear how an agreement forest could capture edges that first get added and then
removed again (or vice versa), as this seems to be beyond embeddings of an agreement forest
into $N$ and $N'$, respectively. 
Furthermore, \cref{fig:SNPR:doublePrune} shows two networks for which every shortest SNPR-sequences
prunes at least one edge twice. A similar problem exists for the subtree prune and regraft
operation on unrooted phylogenetic trees for which a characterisation in terms of agreement forests
appears to be problematic as well~\cite{WM15}.

\pdfbookmark[1]{Acknowledgments}{Acknowledgments}
\subsection*{Acknowledgements}
We thank the anonymous reviewers for their helpful comments.

\providecommand{\etalchar}[1]{$^{#1}$}

\end{document}